\documentclass[a4paper]{amsart}
\usepackage{amsfonts,amsmath,amsthm, amssymb,enumerate,latexsym,url}

\newtheorem{theorem}{Theorem}[section]
\newtheorem{lemma}[theorem]{Lemma}

\newtheorem{question}[theorem]{Question}
\newtheorem{problem}[theorem]{Problem}
\newtheorem{corollary}[theorem]{Corollary}

\newtheorem{prop}[theorem]{Proposition}

\theoremstyle{definition}
\newtheorem{remark}[theorem]{Remark}

\newtheorem{construction}[theorem]{Construction}
\newtheorem{questions}[theorem]{Questions}

\newtheorem{ttable}[theorem]{Table}

\newcommand{\defeq}{{\buildrel{\rm def}\over{\ =\ }}}

\newcommand{\Z}{\mathbb{Z}}

\newcommand{\R}{\mathbb{R}}
\newcommand{\Rn}{\mathbb{R}^n}

\DeclareMathOperator{\spann}{span}
\DeclareMathOperator{\proj}{proj}

\DeclareMathOperator{\hdim}{dim_{\mathsf{H}}}
\DeclareMathOperator{\pdim}{dim_{\mathsf{P}}}
\DeclareMathOperator{\bdim}{dim_{\mathsf{B}}}
\DeclareMathOperator{\ubdim}{{\overline{dim}_{\mathsf{B}}}}

\DeclareMathOperator{\lbdim}{{\underline{dim}_{\mathsf{B}}}}

\newcommand{\cal}{\mathcal}

\def\csillag{({\raise -3pt \hbox{*}}) }
\def\csillagv{({\raise -3pt \hbox{*}})}
\def\ketcsillag{({\raise -3pt \hbox{**}}) }
\def\haromcsillag{({\raise -3pt \hbox{***}}) }
\def\ketcsillagv{({\raise -3pt \hbox{**}})}
\def\haromcsillagv{({\raise -3pt \hbox{***}})}

\def\normp#1{\Arrowvert #1 \Arrowvert_p}

\newcommand{\al}{\alpha}

\newcommand{\de}{\delta}
\newcommand{\eps}{\varepsilon}

\newcommand{\ka}{\kappa}

\newcommand{\su}{\subset}
\newcommand{\sm}{\setminus}

\newcommand{\ds}{\displaystyle}

\begin{document}

\title{Small union with large set of centers}

\author{Tam\'as Keleti}

\date{}



\address
{Institute of Mathematics, E\"otv\"os Lor\'and University, 
P\'az\-m\'any P\'e\-ter s\'et\'any 1/c, H-1117 Budapest, Hungary}

\email{tamas.keleti@gmail.com}

\urladdr{http://www.cs.elte.hu/analysis/keleti}

\thanks{Supported by 
Hungarian Scientific Foundation grant no.~104178.} 

\begin{abstract}
Let $T\su\R^n$ be a fixed set. 
By a scaled copy of $T$ around $x\in\R^n$ we mean a set 
of the form $x+rT$ for some $r>0$.

In this survey paper we study results about the following type of
problems: How small can a set be if it contains a scaled copy of $T$ 
around every point of a set of given size?
We will consider the cases when $T$ is circle or sphere centered at the origin,
Cantor set in $\R$,
the boundary of a square centered at the origin, 
or more generally the $k$-skeleton ($0\le k<n$)
of an $n$-dimensional cube centered at the origin 
or the $k$-skeleton of a more general polytope of $\R^n$.

We also study the case when we allow not only scaled copies but also 
scaled and rotated copies and also the case when we allow only rotated
copies.
\end{abstract}

\maketitle

\section{Introduction}

In this survey paper we study the following type of problems. 
How small can a set be if it contains a scaled copy of a given set around
a large set of points of $\R^n$? 
More precisely:

\begin{problem}
\label{p:general}
Let $T\su\R^n$ be a
fixed set and let $S,B\su\R^n$ be sets such that for every $x\in S$
there exists an $r>0$ such that $x+rT\su B$, in other words $B$ contains
a scaled copy of $T$ around every point of $S$. How small can $B$ be if 
we know the size of $S$?
\end{problem}

In Section~\ref{s:classical} we study the case when $T$ is a circle or sphere 
centered at the origin and we present the classical deep results of 
Stein \cite{St}, Bourgain \cite{Bo}, Marstrand \cite{Ma}, Mitsis \cite{Mi}
and Wolff \cite{Wo97, Wo}. The strongest results of this section 
(Theorem~\ref{t:mitsiswolff} \cite{Mi,Wo} and its corollary 
Theorem~\ref{t:mathedimension}) state that if the Hausdorff dimension
of $S$ is $\hdim S>1$ then $B$ must have positive Lebesgue measure and if 
$\hdim S\le 1$ then $\hdim B\ge \hdim S+n-1$.   

In Section~\ref{s:Cantor} we study the case when $n=1$ and $T$ is a Cantor
set with $0\not\in T$. In this case there are four results: 
{\L}aba and Pramanik \cite{LP} constructed Cantor sets $T\su[1,2]$ for which 
the Lebesgue measure of $B$ must be positive whenever $S$ has positive 
Lebesgue measure;
M\'ath\'e \cite{AM} constructed Cantor sets $T$ for which this is false;
Hochman \cite{Ho} proved that if $\hdim S>0$ then 
$\hdim B>\hdim C+\delta$ for any porous Cantor set $T$, 
where $\de>0$ depends only on $\hdim T$ and $\hdim S$;
and M\'ath\'e noticed that it is a consequence of a recent projection 
theorem of 
Bourgain~\cite{Bo10}  that  $\hdim C>0$ implies
$\hdim B\ge \frac{\hdim S}2$.

In Section~\ref{s:squares} we present the results of Nagy, Shmerkin and the
author \cite{KNS} about the case when $n=2$ and $T$ is 
the boundary or the set of vertices 
of a fixed axis parallel square centered at the origin.
It turns out that in these cases
$B$ can be much smaller than $S$. 
If $S=\R^2$ 
and $T$ is the boundary of the square then
the minimal Hausdorff dimension 
of $B$ is $1$ (Proposition~\ref{p:Hausdorff1}), 
the minimal upper box, lower box and packing dimensions of $B$ are all $7/4$ 
(Theorem~\ref{t:7/4}).
If $S=\R^2$ and $T$ is the set of vertices of the square then the minimal Hausdorff dimension of $B$ is also $1$
and the minimal upper box, lower box and packing dimensions are $3/2$ 
(Theorems~\ref{t:3/4est} and \ref{t:3/4constr}).
For general $S$ when the dimension $s$ of $S$ is given, where dimension can be
Hausdorff, lower box, upper box and packing dimension, 
the smallest possible dimension of
$B$ is determined as a function of $s$ both in the boundary and in the
vertices cases, see Table~\ref{t:square}.
It is remarkable that we get three different functions for different
dimensions in the boundary case.

In Section~\ref{s:cubes} we study the case when $T$ is the $k$-skeleton 
$(0\le k<n)$ of a fixed axis parallel $n$-dimensional cube 
centered at the origin or the
$k$-skeleton of a more general polytope. 
For the $k$-skeleton of axis parallel cubes Thornton \cite{Th} generalized
the above mentioned two-dimensional results for packing and box dimensions
(Theorem~\ref{t:Thornton}), found the estimate 
$\hdim B\ge \max(k,\hdim S -1)$
for Hausdorff dimension
and posed the conjecture the this estimate is sharp. 
This conjecture was proved by Chang, Cs\"ornyei, H\'era and the author 
\cite{CCsHK} not only for cubes but for more general polytopes 
(Theorem~\ref{t:conj}).
To obtain this result first
the smallest Hausdorff dimension of $B$ was determined for any
fixed compact $S$, in other words, instead of $\hdim S$ we fixed $S$ itself
(Theorem~\ref{t:fixedS}).
This was done by showing that for any fixed compact set $S$ if we take
a suitable $B$ in a (Baire category sense) typical way then $\hdim B$ is 
minimal, see the explanation after Theorem~\ref{t:fixedS} for more details.

So far we allowed only scaled copies of $T$. 
In Section~\ref{s:rotated} first we study the modification of
Problem~\ref{p:general} when scaled and rotated copies of $T$ are also 
allowed.
Obviously, allowing rotations can only decrease the minimal dimension of the set $B$. When $T$ is the $k$-skeleton of a cube, it is hard to imagine that rotated copies of $T$ can overlap ``more'' than non-rotated axis-parallel ones. Therefore one might think that allowing rotations of $T$ cannot yield a smaller set $B$. However, this turns out to be false. In fact, 
for any $0\le k<n$ there exists a closed set $B$ of Hausdorff dimension
$k$ that contains the $k$-skeleton of a (rotated) cube centered at every point
of $\Rn$ (Theorem~\ref{t:rotatedcubes}). 
We also study the the modification of
Problem~\ref{p:general} when instead of scaled copies we use only rotated 
copies of $T$ (of fixed size). 
It turns out (Theorem~\ref{t:unitconstruction}) that even 
by using only rotated copies
we can get smaller union than 
using only scaling:
for any $0\le k<n$ there exists a Borel set $B$ of Hausdorff dimension
$k+1$ that contains the $k$-skeleton of a (rotated) \emph{unit} cube  
centered at every point of $\Rn$.
It turns out that this $k+1$ is sharp, $B$ must have Hausdorff dimension
at least $k+1$ (Corollary~\ref{c:atleastkplus1}).
Most of the results of this section are very recent results of 
Chang, Cs\"ornyei, H\'era and the author \cite{CCsHK},
the last mentioned lower estimate for $\hdim B$ is a very recent
result of H\'era, M\'ath\'e and the author \cite{HKM}.

\section{Classical 
results about circle and sphere packing}
\label{s:classical}

How small can a set be if it contains a scaled copy of a given set around
a large set of points of $\R^n$? The first result 
of this type seems to be 
following result of E. Stein.

\begin{theorem}{(Stein (1976) \cite{St})} 
\label{t:stein}
{Let $S\subset\R^n$ ($n\ge 3$) be a set of 
positive Lebesgue measure.
If $B\subset\R^n$ contains a sphere centered at every point of 
$S$, then $B$ has positive Lebesgue measure.}
\end{theorem}

In fact, what Stein proved is a stronger result about the spherical
maximal operator. Let 
\begin{equation}
\label{max}
M_{S^{n-1}}f(x)=\sup_{r>0}\int_{S^{n-1}} |f(x+ry)| d\sigma^{n-1}y \qquad (x\in\R^n),
\end{equation}
where $\sigma^{n-1}$ is the normalized surface measure on $S^{n-1}$.
In other words $M_{S^{n-1}}f(x)$ is the maximal spherical average of $|f|$ around $x$.

\begin{theorem}{(Stein's spherical maximal theorem (1976) \cite{St})} 
\label{t:spherical}
{The spherical maximal operator \eqref{max} is bounded from $L^p(\R^n)$ to $L^p(\R^n)$
if $n\ge 3$ and $p>n/(n-1)$.}
\end{theorem}

Theorem~\ref{t:stein} follows easily from Theorem~\ref{t:spherical}:
we claim that if $B\su\Rn$ contains a 
sphere centered at every point of $S$ and $B$ has Lebesgue measure
zero then so has $S$. 
Indeed, if $f$ is the indicator function of $B$ then by 
Theorem~\ref{t:spherical}, $\normp{M_{S^{n-1}}f}=0$. Since, by definition,
$M_{S^{n-1}}f=1$ on $S$ we obtain that $S$ has measure zero.

Note that this argument also shows that it is enough to assume that around
each point of $S$ there is a sphere that intersects $B$ in a set of
positive $(n-1)$-dimensional measure.

For about ten years it was open if the above theorems hold for $n=2$.
They were settled by Marstrand and Bourgain, independently.

\begin{theorem}{(Bourgain (1986) \cite{Bo}, Marstrand (1987) \cite{Mar})} 
\label{t:bourgainmarstrand}
{Let $S\subset\R^2$ be a set of 
positive Lebesgue measure.
If $B\subset\R^2$ contains a circle centered at every point of 
$S$, then $B$ has positive Lebesgue measure.}
\end{theorem}

\begin{theorem}{(Bourgain's circular maximal theorem (1986) \cite{Bo})} 
\label{t:circular}
{The circular maximal operator $M_{S^{1}}$ defined by \eqref{max} is bounded from $L^p(\R^2)$ to 
$L^p(\R^2)$ for $p>2$.}
\end{theorem}

In fact, Bourgain 
proved his results not only for circles but for any smooth curve with
non-vanishing curvature.

In Theorems \ref{t:stein} and \ref{t:bourgainmarstrand} the set $B$ contains an 
$n$-dimensional family of $n-1$-dimensional surfaces (without big intersection), so it is 
not surprising that $B$ has positive Lebesgue measure. We can also expect that 
smaller set of centers should also guarantee that $B$ has positive measure,
perhaps a set of centers with Hausdorff dimension larger than $1$ is large enough. 
This is indeed the case as it was shown by T. Mitsis and T. Wolff:

\begin{theorem}(Mitsis (1999) \cite{Mi} for $n\ge 3$, Wolff (2000) \cite{Wo} for
$n=2$)
\label{t:mitsiswolff}
Let $S\subset\R^n$ ($n\ge 2$) 
be a set of Hausdorff dimension $>1$. 
If $B\subset\R^n$ contains a sphere centered at every point of $S$, then $B$ has positive Lebesgue measure.
\end{theorem}

In fact, Wolff proved the following more general statement.

\begin{theorem}(Wolff (2000) \cite{Wo})
Let $D\su\R^2\times(0,\infty)$ be a set of Hausdorff dimension larger than $1$
and suppose that $B\su \R^2$ is a set such that for any $(x,r)\in D$ 
a positive measure part of the circle $C(x,r)$ is contained in $B$.
Then $B$ must have positive Lebesgue measure.
\end{theorem}

The following example shows that the above results are sharp.

\begin{theorem}(Talagrand (1980) \cite{Ta})
There exists a planar set of Lebesgue measure zero that contains a circle
centered at every point of a straight line.
\end{theorem}

If the Hausdorff dimension of the set of centers is at most $1$ then
we can estimate the Hausdorff dimension of the union of the circles:

\begin{theorem}(Wolff (1997) \cite{Wo97})
\label{t:wolffdimension}
If $S$ is a Borel set in the plane with $\dim S=s\le 1$ and 
$B$ contains a circle centered at every point of $S$ then $\dim B\ge s+1$. 
\end{theorem}

In fact, Theorem~\ref{t:mitsiswolff} implies 
Theorem~\ref{t:wolffdimension} and its following higher dimensional 
generalization. 
The author heard this simple argument from Andr\'as M\'ath\'e but it is quite
likely that Wolff, Mitsis and others were also aware of this.

\begin{theorem}
\label{t:mathedimension}
If $S\su\Rn$ ($n\ge 2$) with $\dim S=s\le 1$ and 
$B$ contains a sphere centered at every point of $S$ then $\dim B\ge s+n-1$. 
\end{theorem}

\begin{proof}
Fix $\eps\in(0,1)$.
First we construct a compact set $D\su\Rn$ 
with
$\hdim D=\bdim D=1-s+\eps$ such that $\hdim (S+D) = 1+\eps$. 
Let $D_1\su\Rn$ be any compact set (for example a self-similar set)
with $\hdim D_1=\bdim D_1=1-s+\eps$. Then $\hdim (S \times D_1)=1+\eps$,
so its orthogonal projection to 
almost every $n$-dimensional subspace of $\R^{2n}$ also has
Hausdorff dimension $1+\eps$. 
This implies that a suitable affine copy of $D_1$ has all the properties 
we wanted for $D$.

Now  $B+D$ contains a sphere centered at every point of $S+D$ and 
$\hdim (S+D) = 1+\eps>1$, so by Theorem~\ref{t:mitsiswolff}, $B+D$ has
positive Lebesgue measure, so $\hdim(B+D)=n$. 
On the other hand
$\hdim(B+D)\le \hdim B + \bdim D= \hdim B + 1-s+\eps$,
which implies $\hdim B \ge s+n-1-\eps$ for any $\eps\in(0,1)$.
\end{proof}

Note that all of the above results are in harmony with 
the general principle that 
an $a$-dimensional
collection of $b$ dimensional sets in $\Rn$ must have positive measure 
if $a+b>n$ and Hausdorff dimension $a+b$ if $a+b\le n$ unless the
sets have large intersections. 

\section{Cantor sets}
\label{s:Cantor}

In this section we study the case when we want the scaled copy of
a Cantor set around a given set of points. 
(By Cantor set we mean nowhere dense nonempty compact set without 
isolated points.)
 
For any finite Borel 
measure $\mu$ on $\R^n$ one can generalize the maximal operator
\eqref{max} as
\begin{equation}
\label{mumax}
M_{\mu}f(x)=\sup_{r>0}\int_{\R^n} |f(x+ry)| d\mu(y) \qquad (x\in\R^n).
\end{equation}

I.~{\L}aba and M.~Pramanik \cite{LP} studied this maximal operator in the 
case when $\mu$ is the 
natural probability measure on a Cantor set $C\su[1,2]$.
They constructed Cantor sets of Hausdorff dimension $s$ for any 
$s\in(\frac23,1]$ such that $M_{\mu}$ of \eqref{mumax} is bounded on $L^p(\R)$
for any $p>\frac{2-s}s$. Similarly as Theorem~\ref{t:spherical} easily implies
Theorem~\ref{t:stein}, this result gives the following.

\begin{theorem}({\L}aba-Pramanik (2011)  \cite{LP})
For any $s\in(\frac23,1]$ there exists a Cantor set $C\su[1,2]$ of 
Hausdorff dimension $s$ with the following property.
If $B$ and $S$ are subsets of $\R$, $S$ has positive Lebesgue measure
and for every $x\in S$ there exists an $r>0$ such that $rC+x\su B$ 
then $B$ must also have 
positive Lebesgue measure.
\end{theorem}

It is natural to ask if the same is true for every Cantor set.
Very recently A.~M\'ath\'e showed that the answer is negative: 

\begin{theorem}(M\'ath\'e \cite{AM})
For any $s\in(0,1)$ there exists a Cantor set $C$  of 
Hausdorff dimension $s$ and a Borel set $B$ of Hausdorff dimension 
$\frac{s+1}2$ such that 
for every $x\in \R$ there exists an $r>0$ such that $rC+x\su B$.
\end{theorem}

The idea of the proof of the above result of M\'ath\'e comes from a 
discrete grid
construction of Gy.~Elekes \cite{El} that shows the sharpness
of the Szemer\'edi-Trotter theorem \cite{ST}. 

By the above results, for some Cantor sets $C$ there is a set $B\su\R$ of measure zero
which contains scaled copies of $C$ around every $x\in\R$ and for
some other sets $C$ there is no such $B$.
For given Cantor sets $C$ with $0\not\in C$ it seems to be hard to
decide whether 
such a construction is possible or not. 
For example, we do not know this
if $C$ is the classical middle-$\frac13$ Cantor set translated
so that it is symmetric around the origin. 

One can also ask if at least we can guarantee that the union of scaled
copies of a Cantor set $C$ around every point of $\R$ has Hausdorff
dimension strictly larger than the Hausdorff dimension of $C$. 
Very recently M.~Hochman \cite{Ho} gave an affirmative answer to this question
for porous Cantor sets. (Here a set is called
\emph{porous} if there exist $c>0$ and $r_0>0$ such that every
interval of length $r<r_0$ contains an interval of length $cr$ disjoint to
the set.)

\begin{theorem}(Hochman (2016) \cite{Ho})
\label{t:mike}
Let $S\su\R$ be a compact set with $\hdim S>0$,
$C\su\R$ be a porous Cantor set.
If $B\su\R$ contains a  scaled copy $rC+x$ of $C$ for every $x\in S$
then $\hdim B>\hdim C + \de$, where $\de>0$ depends only on $\dim S$
and $\dim C$.
\end{theorem}


Surprisingly, $\hdim B\ge \frac{\hdim S}{2}$ is guaranteed just by assuming $\hdim C>0$. This follows from a recent theorem of Bourgain, as was pointed out by A. M\'ath\'e. 


\begin{theorem}(Bourgain (2010) \cite{Bo10})
\label{t:proj}
For every $0<\al<2$ and $\kappa>0$ there exists $\eta>\al/2$ such that
if $A\su\R^2$ is a set with $\hdim A>\al$ then 
$\hdim \proj_\theta (A)\ge \eta$ for every $\theta\in S^1$ except in an 
exceptional set $E\su S^1$ satisfying $\hdim E\le \ka$.
\end{theorem}

\begin{corollary}
\label{c:1/2}
Let $S$ and $C$ be compact subsets of $\R$ with $\hdim C>0$.
If $B$ is  a Borel subset of $\R$ 
and for every $x\in S$ there exists an $r(x)\in\R$ such that $r(x)C+x\su B$
then $\hdim B\ge \frac{\hdim S}2$.
\end{corollary}

(For a short direct proof of the corollary, see M\'ath\'e \cite{AM}.)


\begin{proof}
We can clearly suppose that $\hdim S>0$.
Let $0<\eps<\frac{\hdim S}{2}$, 
$0<\ka<\hdim C$ and $\al=\hdim S-2\eps$ and let 
$\eta>\al/2=\frac{\hdim S}{2}-\eps$ be the number guaranteed by Theorem~\ref{t:proj}.
Let $A$ be the graph of the function $r(x)$; that is, $A=\{(x,r(x)) : x\in S\}$.
Then $\hdim A\ge \hdim S>\al$.
Let $p_c(x,y)=x+cy$. Note that for every $c\in C$, $p_c(A)\su B$ and 
$p_c(A)$ is a scaled copy of a projection of $A$. 
By the theorem,  
$\hdim p_c (A)\ge \eta$ for every $c\in \R$ except in an 
exceptional set $E$ satisfying $\hdim E\le \ka$.
Since $\hdim E\le \ka<\hdim C$ we can choose $c\in C\sm E$.
Thus $\hdim(B)\ge\hdim(p_c (A)) \ge \eta>\frac{\hdim S}{2}-\eps$.
\end{proof}


We remark that Theorem~\ref{t:mike} and Corollary~\ref{c:1/2} 
are easier for upper box dimension.
In case of Theorem~\ref{t:mike} see the explanation in \cite{Ho}
before the proof the theorem.
In case of Corollary~\ref{c:1/2}, as M\'ath\'e pointed out, 
if we consider upper box dimension instead of
Hausdorff dimension then Corollary~\ref{c:1/2}
would follow very easily even if we assume only 
that $C$ has at least two points. 
Indeed, if $B$ contains a scaled copy $r\{p,q\}+x$ of a set $\{p,q\}$ 
($p\neq q$) for
every $x\in S$ then, as one can easily check, we have 
$\frac{p}{p-q}B-\frac{q}{p-q}B\supset S$, which implies
that $\ubdim B\ge \frac{\ubdim S}2$.


\section{Squares}
\label{s:squares}

The results and arguments of this section are 
due to D.~T.~Nagy, P.~Shmerkin and the author \cite{KNS}.
 
First we consider Problem~\ref{p:general} in the case when $n=2$ and $T$ is
the boundary of a fixed axis-parallel square 
centered at the origin, say $T$ is the boundary of
$[-1,1]\times[-1,1]$. That is, now we study the following problem.

\begin{question}[The question for square boundaries]
\label{q:squareboundary}
Let $T$ be the boundary of the square 
$[-1,1]\times[-1,1]$.
Suppose that we have sets $B, S\su \R^2$ such that
$B$ contains the boundary of an axis parallel 
square centered at every point
of S; that is,
$$
(*)\ (\forall x\in S)\ (\exists r>0)\ x+rT\su B.
$$
How big must $B$ be (or how small can it be) if
the size of $S$ is given? 
\end{question}

It is not hard to see that we can suppose that
$B$ is of the form $B=(A\times \R)\cup (\R\times A)$.
In this case it is easy to see that condition (*)
about $B\su\R^2$ is equivalent to the 
following condition about $A\su\R$:

$$
{\ketcsillag}\quad
{(\forall (x,y)\in S) \ (\exists r\neq 0) \ 
x-r, x+r, y-r, y+r \in A.}
$$

Therefore we obtained the following one-dimensional problem:

\begin{question}[$1$-dimensional problem]
\label{q:1Dproblem}
How small can $A$ be with the following property
if the size of $S$ is given?
$$
{\ketcsillag}\qquad
{(\forall (x,y)\in S) \ (\exists r\neq 0) \ 
x-r, x+r, y-r, y+r \in A.}
$$
\end{question}

It is very easy to check that {\ketcsillag} is equivalent to
$$
{\haromcsillag}\quad
{(A-x) \cap (x-A) \cap (A-y) \cap (y-A) \not\su \{0\}
\quad\text{for all }(x,y)\in S}.
$$

Note that by the Baire category theorem {\haromcsillag} holds even for $S=\R^2$
for any dense $G_\delta$ set $A\su\R$. 
Since there exist dense $G_\delta$ sets of Hausdorff dimension zero,
this means that a set $A$ that satisfies {\haromcsillag} 
(or equivalently {\ketcsillag}) can even have
Hausdorff dimension zero even if $S=\R^2$. 
Thus going back to Question~\ref{q:squareboundary}
we get the following.

\begin{prop}[\cite{KNS}]
\label{p:Hausdorff1}
There exists a set $B\subset\R^2$ that contains a square centered at every 
point of $\R^2$ and has {Hausdorff dimension $1$} (same as a single square!).
\end{prop}

\begin{proof}
Let $A\subset \R$ be a dense $G_\delta$ set of Hausdorff dimension $0$. 
Then $B=(A\times\R)\cup(\R\times A)$ is good (by our previous observations).
\end{proof}

\begin{remark}
A construction of Davies, Marstrand and Taylor \cite{DMT} 
shows that $A$ (and therefore $B$) can be taken to be compact 
for bounded $S$.
\end{remark}

The set $B$ in Proposition~\ref{p:Hausdorff1} is ``small'' from the point of 
view of Hausdorff dimension, but it is ``large'' from the point of view of 
topology and other fractal dimensions. 

The following results show that Question~\ref{q:squareboundary}
becomes more interesting when one considers other notions of 
fractal dimension and/or discrete analogs. 
It turns out that, in some sense, 
Hausdorff dimension is the ``wrong'' dimension for this problem.

First we consider the $S=\R^2$ case.




\begin{theorem}[\cite{KNS}]
\label{t:7/4}
(i) {If $B\su\R^2$ is a set that contains an axis parallel 
square boundary centered at every point of  $[0,1]^2$ then its
lower box, upper box and packing dimension is at least  {$7/4$}.}
  
(ii) {There exists a compact set $B$ that contains a 
square boundary centered at every point of  $[0,1]^2$ such that 
$B$ has lower box, upper box and packing dimension {$7/4$}.}
\end{theorem}

If we let $S$ be a compact subset of $\R^2$ of given dimension 
then even box and packing
dimensions give completely different results:

\begin{theorem}(\cite{KNS})

\noindent
(i) Let $B\su\R^2$ be a set that contains an axis parallel square boundary
centered at every point of  $S\su\R^2$. Then:

(a) If {$\dim=\ubdim$ or $\lbdim$} then 
$\ds\dim B \ge \max \left(1,\frac{7}{8}\dim S\right)$.

(b)
{$\ds\pdim B \ge 1 + \frac{3}{8}\pdim S$}.

\noindent
(ii) The above results are sharp: 
for each $s\in[0,2]$ 

(a') there exist compact sets $S,B$ as above such that 
$$\bdim S = s \qquad \textrm{and} \qquad
\ds\bdim B = \max \left(1,\frac{7}{8}s\right).$$

(b') there exist compact sets $S,B$ as above such that 
$$\pdim S = s \qquad \textrm{and}\qquad \ds\pdim B = 1 + \frac{3}{8}s.$$
\end{theorem}




To get the box and packing dimension estimates
the analogous discrete results were needed.
To get the analogous discrete results the  
following additive number theory lemma was needed.

\begin{lemma}(1-Dimensional Main Lemma \cite{KNS})
\label{l:1Dmain}
If $A\su\R$ and $S\su\R^2$ are finite sets
such that
$(\forall (x,y)\in S)\ (\exists r>0)\  
x-r, x+r, y-r, y+r\in A$ then 
$$
|A| \ge \frac{1}{\sqrt 2}|S|^{\frac{3}{8}}.
$$ 
\end{lemma}

After seeing the connection between the $2$-dimensional problem 
Question~\ref{q:squareboundary} and the $1$-dimensional
problem Question~\ref{q:1Dproblem}, it 
is not surprising that we need this $1$-dimensional lemma to get
the $2$-dimensional discrete result. 
The remarkable fact is that
to get the $1$-dimensional lemma it helps
to go back to $2$-dimension and to prove 
the following.

\begin{lemma}(2-Dimensional Main Lemma \cite{KNS})
\label{l:2Dmain}
If $B, S\su\R^2$ are finite sets
such that
$(\forall (x,y)\in S)\ (\exists r>0)\  
(x\pm r, y\pm r) \in B$ then 
$$
|B| \ge \left(\frac{|S|}{2}\right)^{\frac{3}{4}}. 
$$
\end{lemma}

All discrete, box dimension and packing dimension estimates in \cite{KNS}
are based on this combinatorial geometry lemma. 
Its proof is a rather short double counting argument but it is highly 
nontrivial. 

Although this lemma was motivated by the square boundary problem 
Question~\ref{q:squareboundary},
as a spinoff it can be directly applied to the problems
when we want not the boundary but only the 4 vertices of axis-parallel 
squares. This leads to the following result:

\begin{theorem}{(\cite{KNS})}
\label{t:3/4est}
If $B, S\su\R^2$ are such that $B$ contains the four 
{vertices} of an axis-parallel square centered at
every point of $S$ then 
$$
\ds\dim B \ge \frac{3}{4} \dim S.
$$
if $\dim$ is one of 
{$\pdim$, $\lbdim$ or $\ubdim$}.
\end{theorem}

This is sharp:

\begin{theorem}{(\cite{KNS})}
\label{t:3/4constr}
For each $s\in[0,2]$ there are compact sets $S$ and $B$ as above such that
$\pdim S=\bdim S=s$ and $\pdim B=\bdim B=\frac{3}{4}s$.
When $s=2$ we can choose $S=[0,1]^2$.
\end{theorem}

All the above mentioned results that show the sharpness of the discrete,
box dimension or packing dimension estimates are based on a
simple discrete construction. This construction was found
independently by four students of the E\"otv\"os Lor\'and University:
B.~Bodor, A.~M\'esz\'aros, D.~T.~Nagy and I. Tomon at the
Mikl\'os Schweitzer Mathematical Competition in 2012, where
P. Shmerkin and the author 
posed the
following simple but still highly nontrivial discrete version 
of Question~\ref{q:1Dproblem}.

\begin{question}(Mikl\'os Schweitzer Mathematical Competition 2012)
Call $A\su\Z_n(=\Z / n\Z)$ \emph{rich}, if for every 
$x,y\in\Z_n$ there exists $r\in\Z_n$ for which $x-r,x+r,y-r,y+r\in A$.
What are those $\al$ values for which there exist constants $C_{\al}>0$
such that for any odd $n$ every rich set $A\su\Z_n$ has 
cardinality at least $C_{\al} n^{\al}$?
\end{question}

The answer is $\al\le 3/4$ and the part that these $\al$ values are
good is essentially a special case of the $1$-dimensional main lemma 
(Lemma~\ref{l:1Dmain}) and it can be proved by a much simpler 
but still nontrivial double counting argument. The part that no $\al>3/4$
can be good follows from the following simple construction,
on which the above mentioned other constructions are also based on.

\begin{construction}(Bodor-M\'esz\'aros-Nagy-Tomon, 2012)
\label{students}
Fix a positive integer $k$, let $n=k^4$ and let $A\su\Z_n$ 
consist of those four digit numbers in 
base $k$ that have at least one zero digit. 
Then clearly $A$ has cardinality $O(k^3)=O(n^{3/4})$ and it is easy
to check that $A$ is ``rich'': by choosing the last digit of $r$
we can guarantee that the last digit of $x-r$ is zero, and so on,
each digit of $r$ guarantees one zero digit in $x-r,x+r,y-r,y+r$.  
\end{construction}

The problem when we want only the four vertices of an axis-parallel
square centered at every point of a set of given size is also interesting for 
Hausdorff dimension. So now we study the following:

\begin{question}(Square vertices and Hausdorff dimension)
Let $B\su\R^2$ be a set that contains the four vertices of an axis parallel square
centered at every point of  $S\su\R^2$.
How small can the {Hausdorff dimension} of $B$ be if the 
Hausdorff dimension of $S$ is given?
\end{question}

By projecting to the $x=y$ line we get $\proj B \supset \proj S$, hence clearly
we have
$$
{\hdim B \ge \hdim (\proj B) \ge
\hdim(\proj S) \ge \max(\hdim S -1, 0)}.
$$

Somewhat surprisingly this is sharp:

\begin{theorem} (\cite{KNS})
For any $s\in[0,2]$ there are compact sets $S$
and $B$ as above such that
$\hdim S=s$ and $\hdim B=\max(s-1,0)$.
\end{theorem}


The proof of this result in \cite{KNS} is rather involved.
In the next section we describe a more general result that
can be proved by a simpler method.

To finish this section we summarize the results about squares:

\begin{ttable}(\cite{KNS})
\label{t:square}
If $B\su\R^2$ contains the vertices/boundary of 
axis-parallel squares centered at every point of $S\su\R^2$ of
dimension $s$ (for some dimension) then the best lower bound for the 
dimension (for the same dimension) of $B$ is: 

\begin{center}
\begin{tabular}{|c|c|c|}
 \hline
 dimension & vertices & boundary \\
 \hline
 
 $\pdim$ & $\frac{3}{4}s$ & $1+\frac{3}{8}s$ \\
 $\ubdim$ & $\frac{3}{4}s$ & $\max(1, \frac{7}{8}s)$ \\
 $\lbdim$ & $\frac{3}{4}s$ & $\max(1, \frac{7}{8}s)$ \\
 $\hdim$ & $\max(s-1,0)$ & $1$ \\
 \hline
\end{tabular}
\end{center}
\end{ttable}

\medskip

\section{The $k$-skeletons of cubes and other polytopes of $\Rn$}
\label{s:cubes}

Note that the boundary and the set of vertices
of the square can be considered as
the $1$-dimensional skeleton and the $0$-dimensional  
skeleton of the $2$-dimensional cube. 
So the results of Table~\ref{t:square} can be considered as the 
answers to the $n=2$ special case of the following more general problem.

\begin{question}{($k$-skeleton of $n$-cubes)}
Let $0\le k <n$ and $B\su\R^n$ be a set that contains 
the $k$-dimensional skeleton
of an axis parallel $n$-dimensional cube
centered at every point of  $S\su\R^n$.
How small can the Hausdorff/packing/box dimension of $B$ be if the 
appropriate dimension of $S$ is given?
\end{question}

For box and packing dimension R. Thornton \cite{Th} answered this 
question:

\begin{theorem}(Thornton \cite{Th})
\label{t:Thornton}

(i) For any $0\leq k<n$ and any sets $B,S\su \R^n$ such that 
$B$ contains the $k$-skeleton of an axis parallel 
$n$-cube centered at every point in $S$ we have
$$
(a)\qquad \pdim B\geq k+\frac{(n-k)(2n-1)}{2n^2}\pdim S,
$$
and if $\dim$ denotes upper or lower box dimension then
$$
(b)\qquad \dim B\geq 
\max\left\{k,\left(1-\frac{n-k}{2n^2}\right) \dim S\right\}.
$$

(ii)
Given any $0\leq k<n$, $s\in[0,n]$, there are compact sets
$B,S, B', S'\su\R^n$ where $\pdim(S)=\bdim(S')=s$, $B$ and $B'$ 
contain the $k$-skeleton of an axis parallel $n$-cube centered at every point 
in $S$ and $S'$ respectively, and
$$
(a')\qquad \pdim B=k+\frac{(n-k)(2n-1)}{2n^2}s \quad \textrm{and}
$$
$$
(b')\qquad \bdim B'=\max\left\{k, \left(1-\frac{n-k}{2n^2}\right)s\right\}.
$$
\end{theorem}

The proof of part (ii) is again based on Construction~\ref{students}. 
The proof of part (i) is based on the following generalization of the $1$-dimensional Main Lemma (Lemma \ref{l:1Dmain}) 
and the $2$-dimensional Main Lemma (Lemma \ref{l:2Dmain}).

\begin{lemma}($n,l$-dimensional Main Lemma, Thornton \cite{Th})
\label{l:nlmain}
If $l\le n$ are positive integers, $A\su\R^l$, $S\su\R^n$ and
$$
\forall x\in S \ \exists r>0 \ \forall 1\le i_1<\ldots<i_l\le n\ :\ 
(x_{i_1}\pm r,\ldots,x_{i_l}\pm r)\in A
$$
then 
$$
|A| \ge C_{n,l} |S|^{l(2n-1)/(2n^2)}
$$
for some $C_{n,l}>0$ that depends only on $n$ and $l$.
\end{lemma}

To prove the above lemma,
first the $2$-dimensional Main Lemma (Lemma~\ref{l:2Dmain}) had to be generalized to
$n$-dimension, and then the key observation was that 
Lov\'asz's following corollary of the Katona-Kruskal theorem \cite{Ka,Kr} can be applied:

\begin{theorem}(Katona-Kruskal-Lov\'asz \cite{Lo})
Let $k<n$ be positive integers, $X$ be a finite set of $n$-element sets and let $Y$ be the $k$ element subsets of the sets of $X$.
If $x>0$ are chosen such that $\binom{x}{n}=|X|$, where 
$\binom{x}{n}=x\cdot (x-1) \cdot \ldots \cdot (x-n+1) / n!$,
then $|Y|\ge \binom{x}{k}$.
\end{theorem}

The case of Hausdorff dimension is again completely different. 
Thornton \cite{Th} noticed that a simple projection argument again gives a bound:

If $B$ contains a $k$-skeleton of an axis-parallel $n$-cube centered at every point of
$S$ with $\hdim S=s$ then, denoting the orthogonal projection to the hyperplane with normal vector
$(1,\ldots,1)$ by $\proj$, we have $\proj B\supset \proj S$, so
again $\hdim B\ge \hdim(\proj B) \ge \hdim(\proj S) \ge \hdim S -1=s-1$. 
Since clearly $\hdim B \ge k$, we obtained the following:

\begin{prop} (Thornton \cite{Th})
\label{p:Thorntonestimate}
If $B$ contains a $k$-skeleton of an axis-parallel $n$-cube centered at every point 
of $S$ with $\hdim S=s$ then $\hdim B \ge \max\{k,s-1\}$.
\end{prop}

After proving that this is sharp in some special cases, 
Thornton \cite{Th} posed the conjecture that this is always sharp; that is,
for any integers $0\le k <n$ and any $s\in[0,n]$ there exist compact sets $S$ and $B$
with the above property such that $\hdim B=\max\{k,s-1\}$.

This conjecture was proved by A.~Chang, M.~Cs\"ornyei, K.~H\'era and the author 
\cite{CCsHK}, even for more 
general polytopes:

\begin{theorem}(\cite{CCsHK})
\label{t:conj}
Let $0\le k<n$ and $T$ be the $k$-skeleton of an $n$-dimensional polytope
such that $0$ is not contained in any of the $k$-dimensional affine subspaces
defined by $T$. 
Then for every $s\in[0,n]$ there exist compact sets $S$ and $B$ such that
$$
(\spadesuit)\quad
(\forall x\in S)\ (\exists r>0)\ x+rT\su B,
$$
$\hdim S=s$ and {$\hdim B=\max(s-1,k)$}.
\end{theorem}

We will see later that 
the previous simple projection argument works 
also in this more general case, so 
($\spadesuit$) implies {$\hdim B \ge \max(\hdim S -1, k)$}, therefore the above result shows that this is a sharp estimate for these more general $k$-skeletons.
If $0$ is contained in one of the $k$-dimensional affine subspaces defined by $T$
then the projection argument implies $\hdim B \ge \max(\hdim S, k)$ and it is also
proved in \cite{CCsHK} that in this case this is sharp.
So the minimal Hausdorff dimension of a compact $B$ such that ($\spadesuit$) 
holds
for some compact $S$ with $\hdim S=s$ is $\max\{s,k\}$ if 
$0$ is contained in one of the $k$-dimensional affine subspaces defined by $T$
and $\max\{s-1,k\}$ otherwise.

The following table extends Table~\ref{t:square}
with the results about the $k$-skeletons
an $n$-cube.

\begin{ttable}(\cite{CCsHK},\cite{KNS},\cite{Th})
\label{t:cube}
If $0\le k < n$ and $B\su\R^n$ contains a 
{$k$-dimensional skeleton
of an $n$-dimensional axis-parallel cube} centered at every point 
$S\su\R^n$ of
dimension $s$ (for some dimension) then the best lower bound for the 
dimension (for the same dimension) of $B$ is shown in the last column of the
following table.

\begin{center}
\begin{tabular}{|c|c|c|c|}
 \hline
 dimension & square vertices & squre boundary & 
 {$k$-skeleton of an $n$-cube}\\
           &   ($n=2, k=0$)     &  ($n=2, k=1$) & \\
 \hline
 
 $\pdim$ & $\frac{3}{4}s$ & $1+\frac{3}{8}s$ &
 {$\ds k+\frac{(n-k)(2n-1)}{2n^2}s$}\\
 $\ubdim$ & $\frac{3}{4}s$ & $\ds\max\left\{1, \frac{7}{8}s\right\}$ &
 {$\ds \max\left\{k,\left(1-\frac{n-k}{2n^2}\right)s\right\}$}\\
 $\lbdim$ & $\frac{3}{4}s$ & $\ds\max\left\{1, \frac{7}{8}s\right\}$ &
 {$\ds \max\left\{k,\left(1-\frac{n-k}{2n^2}\right)s\right\}$}\\
 $\hdim$ & $\max\{0,s-1\}$ & $1$ &
 $\max\{k,s-1\}$ \\
 \hline
\end{tabular}
\end{center}

\end{ttable}

\medskip

To get the above mentioned constructions of Chang, Cs\"ornyei, H\'era and the 
author \cite{CCsHK}
a slightly different problem was studied first: 
instead of fixing the Hausdorff dimension $s$ of $S$, the set $S$ itself
was fixed.

\begin{problem}(The modified problem)
Let $T$ and $S$ be fixed subsets of $\R^n$. We want to find a (compact) set $B\su\R^n$
with minimal Hausdorff dimension such that $B$ contains a scaled copy of $T$ around
every point of $S$; that is, 
$$(\spadesuit)\ (\forall x\in S)\ (\exists r>0)\ x+rT\su B.$$
\end{problem}

In \cite{CCsHK} 
this problem was answered for a bit more general $T$
but for simplicity here we consider only the case when $T$ is the $k$-skeleton of a polytope $P$.

Again, we can easily get lower estimate for the Hausdorff dimension using appropriate projections:

 \smallskip

Let $F_1,\ldots,F_m$ be the $k$-dimensional faces of the polytope $P$. 
Let $W_i=\spann F_i$, where $\spann$ means the linear span (not the affine span). 
Then for 
every $u\in\proj_{W_i^\perp} S$ the set $B$ contains a similar copy of $F_i$ in $u+W_i$, which gives $\hdim B \ge \hdim (\proj_{W_i^\perp} S) + k$ for every $i$. 
Therefore we get the estimate
$$
\hdim B \ge \max_i  \hdim (\proj_{W_i^\perp} S) + k \defeq d_{T,S}.
$$
Again, this turned out to be sharp:






\begin{theorem}(\cite{CCsHK})
\label{t:fixedS}
Let $0\le k<n$ and $T$ be the $k$-skeleton of an $n$-dimensional polytope
and let $S\su \R^n$ be a compact set. 
Let $F_1,\ldots,F_m$ be the $k$-faces of $T$ (so $T=\cup_i F_i$), let $W_i=\spann F_i$
and $d_{T,S}=\max_i  \hdim (\proj_{W_i^\perp} S) + k$. 
\begin{itemize}

\item[(a)] 
{If $B$ contains a scaled copy of $T$ around every point of $S$ then
$\hdim B \ge d_{T,S}$.}

\item[(b)] {There exists a compact set $B$ with $\hdim B = d_{T,S}$ that  
contains a scaled copy of $T$ around every point of $S$.}

\end{itemize}

\end{theorem}

We saw the very simple proof of part (a) above. 
Part (b) was shown by proving that in some sense a typical compact $B$
with property ($\spadesuit$) has Hausdorff dimension $d_{T,S}$. 
Now we make this more precise.

We have fixed $T$ and $S$. The set $B$ will be of the form
$$
B=\bigcup_{(x,r)\in K} x+rT \defeq \phi(K),
$$
where $K\su S\times [1,2]$ is compact. 

Note that $B=\phi(K)$ contains a scaled copy of $T$ around every point of $S$ if 
{$\proj_1 K=S$}, where $\proj_1$ denotes the projection to the first coordinate.

Therefore we want to find a ``code set'' $K$ from
$$
{\cal K}=\{K\su S\times [1,2]\ : K \textrm{ is compact}, \proj_1 K=S\}
$$
such that $\hdim (\phi(K))=d_{T,S}$.
By considering the Hausdorff metric, 
$\cal K$ is a compact metric space, so Baire category theorem can be applied.
Thus we can say that 
(in the {Baire category sense}) for a 
{typical} $K\in \cal K$
the set $B=\phi(K)$ has Hausdorff dimension  $d_{T,S}$
if there is a dense $G_{\delta}$ subset $\cal G$ of $\cal K$ 
such that for any $K\in\cal G$ we have $\hdim(\phi(K))=d_{T,S}$.

The point is that this way
it is enough to consider only {one $k$-face} of $T$, and for one $k$-face
it is not hard to show the claim.


To get Theorem~\ref{t:conj} about the original problem in which only the
Hausdorff dimension $s$ of $S$ is given we have to apply 
Theorem~\ref{t:fixedS} for a compact 
set $S$ of Hausdorff dimension $s$ such that the dimension drop of $S$ is maximal for
all the projections $\proj_{W_i^\perp}$.
The existence of such an $S$ can be proved using sets with large intersection
properties defined by K. Falconer \cite{Fa94}.

Finally we explain what the above mentioned maximal dimension drops are
for which $S$ has to be constructed,
how we get the dimension $\max\{k,s-1\}$ or $\max\{k,s\}$
as the dimension of $B$, as it was mentioned after Theorem~\ref{t:conj}, and 
why we have two cases depending on the positions of $0$ and $T$.
The point is that the linear span $W_i$ of the $k$-dimensional face $F_i$ is
$k$-dimensional if the affine span of $F_i$ contains $0$ and $k+1$ dimensional
if not. So the maximal dimension drop of the $s$-dimensional $S$ 
when we project to $W_i^\perp$ is $\min\{s,k\}$ 
in the first case and $\min\{s,k+1\}$ in the second case. Therefore,
if we have maximal dimension drop for each $i$ then
$d_{T,S}=\max_i  \hdim (\proj_{W_i^\perp} S) + k$ equals to
$s-\min\{s,k\}+k=\max\{k,s\}$ 
if we are in the first case for at least one $i$ and 
 $s-\min\{s,k+1\}+k=\max\{k,s-1\}$ 
otherwise.


\section{Rotated cubes and polytopes}
\label{s:rotated}

So far we allowed only scaled copies of a cube or a more general polytope 
but it seems to be natural to study the case when we allow rotations as well. 
For example, how small can $B\su\R^n$ be if it contains the $k$-skeleton 
of a (possibly rotated) $n$-cube centered at every point of 
a set $S$ of given size?

When we did not allow rotation and used only axis-parallel cubes,
we saw that $B$ can be much smaller than $S$
and this was possible because there was a huge overlap between 
the $k$-skeletons of different cubes. This seems to indicate that there
is no point using rotated copies, since their intersection is small. 
This intuition turned out to be wrong, at least for Hausdorff dimension,
as the following recent
result of Chang, Cs\"ornyei, 
H\'era and the author \cite{CCsHK} shows.

\begin{theorem}(\cite{CCsHK})
\label{t:rotatedcubes}
For any $0\le k< n$ there exists a closed set $B$ of Hausdorff dimension $k$
that contains the $k$-skeleton of a (rotated) cube centered at every point
of $\Rn$.
\end{theorem}

Recall from the previous section 
that without allowing rotations the dimension of the above 
Borel set cannot be less than $n-1$. 

We do not know if the results of the previous sections about the box and 
packing dimension are also changed if we allow rotations as well. 

Like the constructions of the previous section for Hausdorff dimension,
Theorem~\ref{t:rotatedcubes} was also shown by Baire category argument. 
Again, we make this more precise.
Let $T$ be the $k$-skeleton of a fixed $n$-dimensional
cube centered at the origin. 
We will put a cube centered at every point of an arbitrary fix compact set $S\su\Rn$.

The set $B$ will be of the form
$$
B=\bigcup_{(x,r, \al)\in K} x+r\al(T) \defeq \Psi(K),
$$
where $K\su S\times [1,2]\times SO(n)$ is compact. 
Now $B=\Psi(K)$ contains a scaled and rotated copy of $T$ centered at every point of 
$S$ if 
{$\proj_1 K=S$}, where $\proj_1$ still denotes the projection to the first coordinate.
Therefore now we want to find a code set $K$ from
$$
{\cal K'}=\{K\su S\times [1,2]\times SO(n)\ : K \textrm{ is compact}, \proj_1 K=S\}
$$
such that $\hdim (\Psi(K))=k$.
By considering a natural compact metric on $SO(n)$,
$\cal K'$ is also a compact metric space with the Hausdorff metric, so Baire category theorem can be applied again.
It is proved in \cite{CCsHK} that for a typical $K\in\cal K'$ indeed we have $\hdim (\Psi(K))=k$, which proves Theorem~\ref{t:rotatedcubes}.

The above argument works of course not only for cubes but for more general sets,
among others for the $k$-skeleton of any polytope provided that 
$0$ is not contained in any of the $k$-dimensional affine subspaces defined by the
polytope.

\medskip

In \cite{CCsHK} we also study what happens if we allow rotation 
but do not allow scaling.
For example, how small the Hausdorff dimension of a set can be if it 
contains the $k$-skeleton of 
an $n$-dimensional (rotated) \emph{unit} cube centered at every point of $\Rn$?
This is closely related to the problems when we want a $k$-dimensional affine 
subspace at
distance $1$ from every point, or when we want $k$-skeleton of 
$n$-dimensional (rotated)
cubes of \emph{every} size centered at every point, or when we want 
$k$-planes at \emph{every} positive distance from every point. 
Combining results of \cite{CCsHK} and \cite{HKM} it turns out that
the smallest possible Hausdorff dimension is $k+1$ in all of these problems.

First we present the result about the constructions.

\begin{theorem}\cite{CCsHK}
\label{t:unitconstruction}
For any integers $0\le k<n$ there exist Borel sets $B_1$ and $B_2$ 
of Hausdorff dimension $k+1$ such that for every $x\in\R^n$ 

(i) the set $B_1$ contains a $k$-dimensional affine subspace at every positive distance from $x$,

(ii) the set $B_2$ contains a $k$-skeleton of (rotated) cubes of every size centered at $x$.
\end{theorem}

Part (i) is very easy: $B_1$ can be chosen as the union of any countable dense 
collection of $k+1$-dimensional affine subspaces.
Like the previously mentioned constructions,
part (ii) is also proved by showing that a typical construction has
Hausdorff dimension at most $k+1$ 
and it also works for the $k$-skeleton of any polytope 
provided that $0$ is not contained in any of the $k$-dimensional affine 
subspaces defined by the polytope.

Recall that, by Proposition~\ref{p:Thorntonestimate}, if instead of $k$-skeleton
of rotated cubes of every size
we are allowed to use axis parallel cubes of arbitrary size then 
$\hdim B\ge n-1$. Therefore, somewhat surprisingly,
Theorem~\ref{t:unitconstruction} (ii) shows that in the $k<n-2$ case
we can get smaller union if we take 
$k$-skeletons of rotated cubes of every size centered at every point
than in the case when we take $k$-skeletons of axis parallel cubes
of arbitrary size centered at every point.

The fact that we cannot get sets of Hausdorff dimension less than $k+1$
even if we want only $k$-skeleton of unit cubes around every
point or $k$-dimensional affine subspace at distance $1$ from every point
is based on the following very recent result of H\'era, M\'ath\'e and
the author.

\begin{theorem}(\cite{HKM})\label{t:dimofunion}
Let $0\le k<n$ be integers and let
$A(n,k)$ denote the space of all $k$-dimensional affine subspaces of $\R^n$
and consider any natural metric on $A(n,k)$. 

Suppose that $B\su\R^n$, $E\su A(n,k)$ and every $k$-dimensional 
affine subspace $P\in E$ intersects $B$ in a set of 
positive $k$-dimensional Hausdorff measure. 
Then
\begin{equation}\label{e:estimate}
\hdim B \ge k+\min(\hdim E,1).
\end{equation}
\end{theorem}

For the special case when $k=n-1$ this was proved by 
Falconer and Mattila \cite{FM}. Note that for $0<k<n-1$ 
the right-hand side of (\ref{e:estimate}) cannot be replaced by 
the more natural $k+\min(\hdim E, n-k)$: 
if $B$ is a $k+1$-dimensional affine
subspace and $E$ is the set of all $k$-planes of $B$ then  
$\hdim E=k+1>1$, so $\hdim B=k+1<k+\min(\hdim E, n-k)$.

\begin{corollary}(\cite{CCsHK})\label{c:atleastkplus1}
Let $0\le k<n$ be integers, $d>0$ be fixed and 
$B$ be a subset of $\R^n$ such that for every $x\in\R^n$ there exists
a $k$-dimensional affine subspace $P$ at distance $d$ from $x$ such that
$P$ intersects $B$ in a set of positive $k$-dimensional Hausdorff measure. 
Then $\hdim B \ge k+1$.
\end{corollary}

\begin{proof} (Sketch)
Let $E$ be the set of those  $k$-dimensional 
affine subspaces that intersect $B$ in a set of 
positive $k$-dimensional Hausdorff measure. 
By Theorem~\ref{t:dimofunion} it is enough to prove that $\hdim E\ge 1$.
For each $P\in E$ let $D(P)$ be the set
of those $k$-dimensional affine subspaces that are parallel to $P$ and are at
distance $d$ from $P$ (in the Euclidean distance of $\R^n$), and let $C(P)$
be the union of the $k$-planes of $D(P)$. Note that $C(P)$ is a cylinder, 
and it is exactly
the set of those points of $\R^n$ that are at distance $d$ from $P$. 
Thus, by assumption, $\cup_{P\in E} C(P)=\R^n$. 

One can easily check that
for any $P$ we have $\hdim D(P)=n-k-1$. From this it is not hard to
show that $\hdim(\cup_{P\in E} D(P))\le \hdim E + n-k-1$ and then
$\hdim(\cup_{P\in E} C(P))\le \hdim(\cup_{P\in E} D(P)) +k \le \hdim E + n-1$.
Since $\cup_{P\in E} C(P)=\R^n$ this implies that indeed  $\hdim E\ge 1$.
\end{proof}

Combining Theorem~\ref{t:unitconstruction} and Corollary~\ref{c:atleastkplus1}
we immediately get the following four statements.

\begin{corollary}\label{c:four}
For any integers $0\le k<n$ and any of the following four properties
the minimal Hausdorff dimension of a Borel set $B\su\R^n$ with that property
is $k+1$.

(i) $B$ contains the $k$-skeleton of a (rotated) unit cube centered at every 
point of $\R^n$.

(ii) $B$ contains the $k$-skeleton of (rotated) cubes of every size 
centered at every point of $\R^n$.

(iii) $B$ contains a $k$-dimensional affine subspace at distance $1$ from 
every point of $\R^n$.

(iv) $B$ contains a $k$-dimensional affine subspace at every positive distance from 
every point of $\R^n$.
\end{corollary}

For $k=n-1$ in Corollaries~\ref{c:atleastkplus1} and \ref{c:four} 
it is natural to ask if instead of full
Hausdorff dimension we can also guarantee positive Lebesgue measure.
The following results give negative answers in case of 
Corollary~\ref{c:atleastkplus1} and Corollary~\ref{c:four} (i) and (iii).

\begin{theorem}\cite{CCsHK}
For any $n\ge 2$ there exist Borel sets $B_1, B_2\su\R^n$ of Lebesgue 
measure zero such that for every $x\in \R^n$ 

(i) the set $B_1$ contains an $n-1$-dimensional hyperplane at distance $1$ from $x$.

(ii) the set $B_2$ contains the boundary of a (rotated) unit cube centered at $x$.
\end{theorem}

These results are also proved by showing that a typical construction has
Lebesgue measure zero 
and the proof of (ii) also works for any polytope 
provided that $0$ is not contained in any of the $n-1$-dimensional affine 
subspaces defined by the polytope.



\section{Open problems}
\label{s:problems}


Here we collect some of those problems that seem to be open. 
Most of these were already mentioned in the previous sections.

\begin{questions}
Let $C$ be the classical middle-third
Cantor set translated by $-\frac12$ to move its center to $0$.
Let $B\su\R$ be a set such that for any $x\in \R$ there exists $r>0$ such that 
$x+rC\su B$. 
How small can $B$ be? How small can $\hdim B$ be? Can $B$ have Lebesgue measure zero?
What can we say about other Cantor sets $C$? 
\end{questions}

Recall from Section~\ref{s:Cantor} that by the 
theorem of {\L}aba and Pramanik \cite{LP} there exist 
Cantor sets for which such a $B$ must have positive Lebesgue measure,
by the result of M\'ath\'e \cite{AM} there exist
Cantor sets for which such a $B$ can have zero measure, 
and by the results of Hochman \cite{Ho} and Bourgain~\cite{Bo10} 
for any porous Cantor set $C$ with $\hdim C>0$
such a $B$
must have Hausdorff dimension strictly larger than $\hdim C$ and at least $1/2$.

\begin{questions}
What can we say about the box and packing dimension in the problems of Section~\ref{s:rotated}? In other words, how the box and packing dimension results about
squares and cubes of Sections~\ref{s:squares} and \ref{s:cubes} are changed if we allow rotated squares and cubes?
For example, is there a compact subset of the plane that contains a (rotated) 
square boundary centered at every point of $[0,1]\times[0,1]$ with box or packing dimension
less than $7/4$ (which is the smallest possible box and packing dimension 
for axis parallel squares by Theorem~\ref{t:7/4})?
\end{questions}

Theorem~\ref{t:rotatedcubes} shows that in some cases
we can get sets with smaller Hausdorff
dimension if we allow rotation but we do not know if smaller box or smaller packing dimension can be also obtained by allowing rotated squares or cubes.

\bigskip

\noindent\emph{Acknowledgment.}
The author is grateful to Andr\'as M\'ath\'e for helpful discussions and
for checking the paper very carefully.

\end{document}